\documentclass{amsart}

\usepackage{amsmath,amssymb,amsthm}
\usepackage{enumitem}
\usepackage[colorlinks=true,linkcolor=blue,citecolor=blue,urlcolor=blue]{hyperref}

\theoremstyle{plain}
\newtheorem{theorem}{Theorem}[section]
\newtheorem{proposition}[theorem]{Proposition}
\newtheorem{lemma}[theorem]{Lemma}
\newtheorem{corollary}[theorem]{Corollary}
\newtheorem{introthm}{Theorem}

\theoremstyle{definition}
\newtheorem{definition}[theorem]{Definition}

\newtheorem{problem}[theorem]{Problem}

\theoremstyle{remark}
\newtheorem{remark}[theorem]{Remark}

\newcommand{\Fp}{\mathbb{F}_p}
\newcommand{\Z}{\mathbb{Z}}
\newcommand{\Aut}{\operatorname{Aut}}

\newcommand{\UAut}{\operatorname{UAut}}

\newcommand{\gr}{\operatorname{gr}}
\newcommand{\ad}{\operatorname{ad}}
\newcommand{\cl}{\operatorname{cl}}
\newcommand{\I}{\mathrm{I}}
\newcommand{\fl}{\mathfrak l}
\newcommand{\fe}{\mathfrak e}
\newcommand{\fv}{\mathfrak v}
\newcommand{\fa}{\mathfrak a}
\newcommand{\fI}{\mathcal I}
\newcommand{\RR}{\mathrm{RR}}
\newcommand{\Lev}{\mathrm{Lev}}
\newcommand{\Vt}{\widetilde V}
\newcommand{\st}{\tilde s}
\newcommand{\Gt}{\widetilde G}

\title[Short Jennings series and the isomorphism problem]{Short Jennings series and the modular isomorphism problem}
\author{Diego Tuzzolo}
\address{Instituto de Matem\'atica e Estat\'istica, Universidade de S\~ao Paulo, Rua do Mat\~ao, 1010, CEP 05508-090, S\~ao Paulo, SP, Brazil}
\email{diegotuzzolo@usp.br}


\begin{document}

\begin{abstract}
Let $p$ be an odd prime and $G$ a finite $p$-group whose Jennings series
has length at most $p$, that is, $D_{p+1}(G)=1$; when $D_p(G)\ne1$ assume
also that $G/D_p(G)$ has nilpotency class at most $p-2$. We prove that $G$
is determined up to isomorphism by its modular group algebra $\Fp G$, and
indeed by the quotient $\Fp G/\I^{\,\nu+1}$, where $\I$ is the augmentation
ideal and $\nu$ the length of the Jennings series; for an arbitrary finite $p$-group $G$, the quotient $\Fp G/\I^{\,p}$
determines $G/D_p(G)$; and that every
automorphism of $\Fp G$ is the product of an automorphism of $G$ and an
automorphism inducing the identity on $\I/\I^2$. In particular the modular
isomorphism problem has a positive answer for all $p$-groups of exponent
$p$ and class less than $p$, and for all $p$-groups of class at most $p-2$
and exponent dividing $p^2$ whose derived subgroup has exponent $p$.
Consequently every counterexample in odd characteristic has class at
least $p-1$, or a derived subgroup of exponent at least $p^2$, or an
element of order $p^3$. The proof combines a descent theorem in the
truncated group algebra of a free group, valid for every prime, with a
relative rigidity theorem for the unipotent automorphisms of truncated
group algebras of groups of exponent $p$, proved through the Lazard
correspondence. The bound $p$ is the range in which the logarithm of the
group algebra exists.
\end{abstract}

\maketitle

\section{Introduction}\label{sec:intro}

Let $p$ be a prime, $K=\Fp$ the field with $p$ elements, and $G$ a finite
$p$-group. The \emph{modular isomorphism problem} asks whether $G$ is
determined up to isomorphism by the $K$-algebra $KG$. Positive answers are
known for abelian groups \cite{Des56}, groups of order at most $p^5$
\cite{Pas65,SS96}, metacyclic groups \cite{Bag88,San96}, groups of class
$2$ with elementary abelian derived subgroup \cite{San89}, groups with
$[G:Z(G)]\le p^2$ \cite{Dre89}, odd $p$-groups with $D_4(G)=1$
\cite{Her07}, quotients $F/R$ of a free group $F$ with $R$ lying between
two consecutive terms of the Jennings series of $F$ \cite{Roh87,Roh90},
and several further classes surveyed in \cite{Mar22}. The problem was
settled negatively for $p=2$ by Garc\'ia-Lucas, Margolis and del R\'io
\cite{GLMdR22}; for odd $p$ it remains open, and groups of exponent $p$
are listed in \cite{Mar22} among the natural test cases.

The augmentation ideal $\I$ of $KG$ is nilpotent, and its powers cut out
in $G$ the \emph{Jennings series} $D_n(G)=G\cap(1+\I^n)$, given by
Jennings' formula
\begin{equation}\label{eq:jennings}
D_n(G)=\prod_{ip^{\,j}\ge n}\mho_j\bigl(\gamma_i(G)\bigr),
\end{equation}
where $\gamma_i$ is the lower central series and
$\mho_j(X)=\langle x^{p^j}:x\in X\rangle$ \cite{Jen41}. Its length is
$\nu(G)=\max\{n:D_n(G)\ne1\}$. The present paper proves that the modular
isomorphism problem has a positive answer whenever this series is short.

\begin{introthm}\label{thmA}
Let $p$ be an odd prime and $G$ a finite $p$-group with $D_{p+1}(G)=1$.
If $D_p(G)\ne1$, assume moreover that $G/D_p(G)$ has nilpotency class at
most $p-2$. Put $\nu=\nu(G)$. Then:
\begin{enumerate}[label=\textup{(\arabic*)}]
\item if $\Gt$ is a finite $p$-group with
$\Fp\Gt/\I(\Fp\Gt)^{\nu+1}\cong\Fp G/\I(\Fp G)^{\nu+1}$, then
$\Gt/D_{\nu+1}(\Gt)\cong G$; in particular $\Fp\Gt\cong\Fp G$ implies
$\Gt\cong G$;
\item every automorphism of $\Fp G$ is the product of an automorphism of
$G$ and an automorphism inducing the identity on $\I/\I^2$.
\end{enumerate}
\end{introthm}

By \eqref{eq:jennings}, $D_p(G)=\gamma_p(G)\,\mho_1(G)$ and
$D_{p+1}(G)=\gamma_{p+1}(G)\,\mho_1(\gamma_2(G))\,\mho_2(G)$. Hence
$D_p(G)=1$ if and only if $G$ has exponent $p$ and class less than $p$,
in which case the additional hypothesis is vacuous, and we obtain:

\begin{corollary}\label{cor:expp}
Let $p$ be odd. Every finite $p$-group of exponent $p$ and nilpotency
class less than $p$ is determined by its modular group algebra.
\end{corollary}

\begin{corollary}\label{cor:expp2}
Let $p$ be odd, and let $G$ be a finite $p$-group of class at most $p-2$
and exponent dividing $p^2$ whose derived subgroup has exponent $p$. Then
$G$ is determined by its modular group algebra.
\end{corollary}

The tower argument behind Theorem~\ref{thmA} also gives a statement with no
hypothesis on $G$ at all.

\begin{corollary}\label{cor:belowceiling}
Let $p$ be odd and let $G$ be any finite $p$-group. If $\Gt$ is a finite
$p$-group with $\Fp\Gt/\I(\Fp\Gt)^{p}\cong\Fp G/\I(\Fp G)^{p}$, then
$\Gt/D_p(\Gt)\cong G/D_p(G)$. In particular, if $\Fp G\cong\Fp\Gt$, then
$G/D_p(G)\cong\Gt/D_p(\Gt)$.
\end{corollary}

\begin{introthm}\label{thmC}
Let $p$ be odd and let $G,\Gt$ be non-isomorphic finite $p$-groups with
$\Fp G\cong\Fp\Gt$. Then at least one of the following holds:
\begin{enumerate}[label=\textup{(\Roman*)}]
\item $\gamma_{p+1}(G)\ne1$; then $\nu(G)\ge p+1$;
\item $\gamma_2(G)$ has exponent at least $p^2$; then $\nu(G)\ge2p$;
\item $G$ has exponent at least $p^3$; then $\nu(G)\ge p^2$;
\item $D_{p+1}(G)=1\ne D_p(G)$ and $G/D_p(G)$ has class exactly $p-1$;
then $\nu(G)=p$, $\cl(G)\ge p-1$, and $\cl(G)=p$ if $G$ has exponent $p$.
\end{enumerate}
Moreover $G/D_p(G)\cong\Gt/D_p(\Gt)$. In particular $G$ has class at least $p-1$, or a derived subgroup of
exponent at least $p^2$, or an element of order $p^3$.
\end{introthm}

\subsection*{Context}
In the language of the Jennings series, Hertweck's theorem \cite{Her07}
says that $\nu(G)\le3$ implies a positive answer for odd $p$, and
Theorem~\ref{thmA} says the same for $\nu(G)\le p$, with the proviso in
part (IV) of Theorem~\ref{thmC} at the ceiling. For $p\ge5$ Theorem~\ref{thmA}
therefore contains the case $D_4(G)=1$ of Hertweck's theorem; for $p=3$
his result is the stronger one. R\"ohl's theorem \cite[Theorem~3.1.2]{Roh90},
announced in \cite{Roh87}, treats $G=F/R$ with $M_{n+1}(F)\subseteq R
\subseteq M_n(F)$ for some $n$, $F$ free and $M_n(F)$ its Jennings series:
its hypothesis is that the relations of $G$ occupy a single layer of the
Jennings filtration, with no restriction on $n$; Theorem~\ref{thmA}
restricts $n\le p$ instead and asks nothing of the shape of the
relations. The smallest groups in Corollary~\ref{cor:expp} outside all
previous results have order $p^6$ (Section~\ref{sec:examples}).

Part (2) of Theorem~\ref{thmA} is the \emph{Levi factorization}
$\Aut(\Fp G)=\Aut(G)\cdot\UAut(\Fp G)$, in the terminology of
\cite{Roh90,HS07}, where $\UAut$ denotes the automorphisms inducing the
identity on $\I/\I^2$ (R\"ohl's \emph{unipotent automorphisms}). It was
known for relatively free groups \cite[Theorem~3.2.1]{Roh90}. It is not
proved here directly for each $G$: it is inherited from $G/D_\nu(G)$, and
what makes the inheritance work is the relative rigidity theorem of
Section~\ref{sec:rigid}.

\subsection*{Two questions from the literature}
The proof answers, in its range, two questions raised explicitly in
print. R\"ohl \cite[Theorem~3.3]{Roh87} proved for \emph{integral} group
rings that if $F/R$ is characterized then so is $F/[F,R]R^p$, by lifting
an isomorphism to the completed augmentation ideal of a free group and
descending; and he observed \cite[p.~863]{Roh87} that the same induction
would work modularly once one could establish a correspondence of ideals
for $\Fp G\cong\Fp H$, or find a suitable substitute. Ideals of modular
group algebras need not correspond, by an example of Sandling recorded in
\cite[\S3.2]{Roh90}. The substitute is Lemma~\ref{lem:transfer} below:
the ideal generated by the top Jennings layer $D_n(G)$ transfers along
every algebra isomorphism \emph{modulo $\I^{n+1}$}, for graded reasons,
and modulo $\I^{n+1}$ turns out to be enough. This is the descent theorem
of Section~\ref{sec:descent}, valid for every prime.

The second question is Hertweck's. His proof that $G/D_4(G)$ is
determined \cite{Her07} modifies a Jennings basis by terms with a factor
$\tfrac12$, and he remarks that this reflects the absence of
mutually inverse maps $\exp$ and $\log$ in characteristic $p$, while
expressing the hope that the method would lead to connections between
group and Lie-algebra objects. Section~\ref{sec:rigid} provides such a
connection: the unipotent automorphisms of the truncated algebra
$\Fp H/\I^{\,n+1}$ do have a logarithm as soon as $n+1\le p$, this
logarithm is a genuine derivation, and it can be projected onto the
Lazard Lie algebra of the universal central extension of $H$, where it
integrates to an automorphism of the group. This is why every use of
the hypotheses ``$p$ odd'' and ``$\nu\le p$'' in the paper is inside
Section~\ref{sec:rigid}.

\subsection*{The mechanism}
Fix a $p$-group $H$. Hertweck and Soriano \cite{HS07} parametrize the
central extensions $1\to A\to G\to H\to1$ with $A\le\Phi(G)$ by subspaces
$U$ of a \emph{relation module} $V$, isomorphism classes of extensions
corresponding to $\Aut(H)$-orbits. Given $G$ with top Jennings layer
$A=D_n(G)$ and base $H=G/A$, an isomorphism $\Fp G/\I^{\,n+1}\cong
\Fp\Gt/\I^{\,n+1}$ lifts to a unipotent automorphism $\Gamma$ of the
truncated group algebra of a free group $F$ carrying the relations of $G$
to those of $\Gt$ (R\"ohl's lifting principle). By
Lemma~\ref{lem:transfer}, $\Gamma$ also carries the relations of the base
to those of the base of $\Gt$, so the bases coincide once the base is
known to be determined by its truncated algebra; after normalization,
$\Gamma$ preserves the relations of $H$ and acts on $V/W_{n+1}$, where
$W_\bullet$ is the depth flag of $V$, carrying the subspace of $G$ to
that of $\Gt$. If this action is realized by an automorphism of $H$, the
two subspaces lie in one $\Aut(H)$-orbit and $G\cong\Gt$
(Theorem~\ref{thm:descent}). In odd characteristic, for $H$ of exponent
$p$ and class at most $p-2$ and $n\le p$, the action \emph{is} realized
by automorphisms of $H$ (Theorems~\ref{thm:trunc} and~\ref{thm:ceiling}),
and an induction along the Jennings series proves Theorem~\ref{thmA}.

Everything is over the prime field. The reason is
Lemma~\ref{lem:PS}, the recognition of the depth flag of $V$ inside the
small group ring, which rests on the identity $\lambda(g-1)\equiv
g^{\ell}-1$ modulo $\I^{n+1}$ for $\lambda\in\Fp$ and $\ell\in\Z$
representing $\lambda$; Margolis and Sakurai \cite{MS25} explain why this
has no analogue over larger fields.

\subsection*{Organization}
Section~\ref{sec:prelim} collects the tools. Section~\ref{sec:descent}
sets up the truncated free frame, proves the transfer lemma and the
descent theorem, and shows how the Levi factorization propagates; it is
valid for every prime. Section~\ref{sec:rigid} proves the relative
rigidity theorem for bases of exponent $p$, first below the ceiling and
then at it. Section~\ref{sec:main} assembles the main theorems.
Section~\ref{sec:examples} discusses examples, the consistency of the
descent with the counterexample at $p=2$, and the scope of the results;
Section~\ref{sec:open} lists open problems.

\subsection*{Conventions}
Throughout, $p$ is a prime and $K=\Fp$; all groups are finite $p$-groups
unless said otherwise. For a group $G$, $\I=\I(KG)$ is the augmentation
ideal, $D_n(G)=G\cap(1+\I^n)$, $L_n(G)=D_n(G)/D_{n+1}(G)$,
$L(G)=\bigoplus_nL_n(G)$ and $\nu(G)=\max\{n:D_n(G)\ne1\}$. Commutators
are $[g,h]=g^{-1}h^{-1}gh$; $\Phi(G)$ is the Frattini subgroup and
$d(G)=\dim_KG/\Phi(G)$. For a free group $F$ we write $M_n(F)$ for the
right-hand side of \eqref{eq:jennings} with $G=F$; by Jennings' theorem
applied to the finite $p$-group $F/R$, one has $D_n(F/R)=M_n(F)R/R$ for
every $R\lhd F$ of $p$-power index.

\section{Preliminaries}\label{sec:prelim}

\subsection{Jennings, Quillen, graded rigidity}
By Jennings' theorem \cite{Jen41}, the Jennings series of a finite
$p$-group $G$ is given by \eqref{eq:jennings}; in particular $D_2(G)=
\Phi(G)$, $\gamma_n(G)\le D_n(G)$, $[D_m,D_n]\le D_{m+n}$ and
$D_n^{\,p}\le D_{np}$, so that $L(G)$ is a graded Lie algebra over $K$
under commutation, with a $p$-map $L_n\to L_{np}$ induced by $p$-th
powers. By Quillen \cite{Qui68}, the graded algebra $\gr(KG)=
\bigoplus_n\I^n/\I^{n+1}$ is the restricted enveloping algebra $u(L(G))$,
the map $L_n(G)\to\I^n/\I^{n+1}$, $gD_{n+1}\mapsto g-1+\I^{n+1}$, being
an embedding; and $L(G)$ is generated in degree one as a restricted Lie
algebra, since by \eqref{eq:jennings} every element of $D_n(G)$ is a
product of $p^j$-th powers of commutators of weight $i$ with $ip^j\ge n$.

\begin{lemma}[graded rigidity]\label{lem:grrigid}
Let $\Theta:KG\to K\Gt$ be a $K$-algebra isomorphism. Then $\gr\Theta$
maps $L_n(G)$ onto $L_n(\Gt)$ for every $n$; in particular
$\dim L_n(G)=\dim L_n(\Gt)$ and $\nu(G)=\nu(\Gt)$.
\end{lemma}

\begin{proof}
$KG$ is local with radical $\I$, so $\Theta(\I)=\widetilde\I$ and
$\gr\Theta$ is an isomorphism of graded algebras. It maps
$L_1(G)=\I/\I^2$ onto $L_1(\Gt)$, hence the restricted Lie subalgebra
generated by $L_1(G)$, which is $L(G)$, onto $L(\Gt)$, preserving degrees.
\end{proof}

\subsection{The Hertweck--Soriano frame}\label{subsec:HS}
Fix a finite $p$-group $H$ with $d=d(H)$ and a minimal presentation
$1\to R\to F\to H\to1$, $F$ free on $x_1,\dots,x_d$, so that $R\le\Phi(F)$.
Following \cite[\S4]{HS07} put
\[
L=[F,R]R^{\,p},\qquad E=F/L,\qquad V=R/L .
\]
Then $V$ is central and elementary abelian in $E$, contained in
$\Phi(E)$, and $1\to V\to E\to H\to1$ is universal among the
\emph{central Frattini extensions} $1\to A\to G\to H\to1$, those with
$A$ elementary abelian, central, and contained in $\Phi(G)$: each is
$E/U\to H$ for a subspace $U\le V$ \cite[Lemma~4.2]{HS07}. For a finite
$p$-group $E$ with a central elementary abelian subgroup $V$, the
\emph{small group ring} is
\[
s(E,V)=KE/\I(KV)\I(KE),
\]
and $\{V-1\}$ denotes the image of $\I(KV)KE$ in it. By
\cite[(3.1)]{HS07}, $v\mapsto v-1$ induces an isomorphism of abelian
groups $V\cong\{V-1\}$; and $\{V-1\}$ is a central ideal with
\begin{equation}\label{eq:trunc0}
\{V-1\}\cdot\I(s)=\I(s)\cdot\{V-1\}=0 ,
\end{equation}
because $(v-1)t(g-1)$ and $(g-1)(v-1)t=(v-1)(g-1)t$ lie in
$\I(KV)\I(KE)$ for $v\in V$ and $g,t\in E$. The \emph{depth flag} of $V$
is $W_n=V\cap D_n(E)$; note $W_2=V$ when $V\le\Phi(E)$.

\begin{lemma}[the depth flag in the small group ring]\label{lem:PS}
Let $E$ be a finite $p$-group and $V$ a central elementary abelian
subgroup, $s=s(E,V)$. Then $\{V-1\}\cap\I(s)^n=\{W_n-1\}$ for every
$n\ge1$.
\end{lemma}

\begin{proof}
Choose a basis $v_1,\dots,v_k$ of $V$ such that, for every $m$, the
$v_i$ lying in $W_m$ form a basis of $W_m$; say $v_i\in W_{w_i}\setminus
W_{w_i+1}$. Since $W_m/W_{m+1}\to L_m(E)$ is injective, the classes of
the $v_i$ with $w_i=m$ are independent in $L_m(E)$, and we may extend
the $v_i$ to elements $g_1,\dots,g_r$ of $E$ such that, for every $m$,
the $g_j$ with $g_j\in D_m(E)\setminus D_{m+1}(E)$ map to a basis of
$L_m(E)$. By Jennings' theorem \cite{Jen41}, \cite[Chapter~3, \S3]{Pas77},
the products $\prod_j(g_j-1)^{e_j}$, $0\le e_j<p$, taken in a fixed
order, form a basis $\mathcal M$ of $KE$, and the products of weight
$\sum_je_jw(g_j)\ge n$ form a basis of $\I(KE)^n$. Let $\mathcal M'
\subseteq\mathcal M$ be the set of products of length $\sum_je_j\ge2$
involving at least one factor $v_i-1$. We claim that $\I(KV)\I(KE)$ is
the span of $\mathcal M'$. Indeed $\I(KV)\I(KE)$ is spanned by the
$(a-1)m$ with $1\ne a\in V$ and $m\in\mathcal M$ of length $\ge1$;
writing $a=\prod v_i^{c_i}$ and expanding, $a-1$ is a combination of
products of the commuting elements $v_i-1$ of length $\ge1$ and exponents
$<p$; since the $v_i-1$ are central and $(v_i-1)^p=v_i^{\,p}-1=0$, each
product $(\text{such a monomial})\cdot m$ is either $0$ or, after
reordering, an element of $\mathcal M'$. Conversely an element of
$\mathcal M'$ is $(v_i-1)m'$ with $m'\in\I(KE)$. This proves the claim.

Now let $v=\prod v_i^{\,c_i}\in V$ with $v-1\in\I(KE)^n+\I(KV)\I(KE)$.
Expanding, $v-1\equiv\sum_ic_i(v_i-1)$ modulo the span of $\mathcal M'$.
Since $\I(KE)^n$ and $\I(KV)\I(KE)$ are spanned by subsets of the basis
$\mathcal M$, and the $v_i-1$ are basis elements of length one, not in
$\mathcal M'$, each $v_i-1$ with $c_i\ne0$ has weight $w_i\ge n$. Hence
$v\in W_n$. The reverse inclusion is clear.
\end{proof}

\subsection*{The action of $\Aut(H)$ on $V$}
Let $\beta\in\Aut(H)$. Choose $\theta\in\operatorname{End}(F)$ with
$\pi\theta=\beta\pi$, where $\pi:F\to H$ is the presentation. Then
$\theta(R)\subseteq R$, and $\theta$ is surjective modulo $\Phi(F)$,
hence surjective and, $F$ being Hopfian, an automorphism; comparing
orders, $\theta(R)=R$, whence $\theta(L)=L$ and $\theta$ induces
$\hat\beta\in\Aut(E)$ with $\hat\beta(V)=V$, inducing $\beta$ on $H$.

\begin{lemma}\label{lem:welldef}
Let $E$ be a finite $p$-group, $V\le\Phi(E)$ central and elementary
abelian, $H=E/V$. Two automorphisms of $E$ preserving $V$ and inducing
the same automorphism of $H$ agree on $V$.
\end{lemma}

\begin{proof}
Their quotient $\gamma$ satisfies $\gamma(x)=xv_x$ with $v_x\in V$
central of order dividing $p$; so $\gamma$ fixes every $[x,y]$ and every
$x^p$, hence $\Phi(E)\supseteq V$ pointwise.
\end{proof}

Thus $\rho(\beta):=\hat\beta|_V$ is well defined, $\rho:\Aut(H)\to GL(V)$
is a homomorphism, and $\rho(\beta)$ preserves every $W_n$ (since
$\hat\beta$ preserves $D_n(E)$). The easy half of the Hertweck--Soriano
correspondence \cite[Theorem~4.4]{HS07} is what we shall use: for
$U\le V$ and $\beta\in\Aut(H)$, $\hat\beta$ induces an isomorphism
$E/U\cong E/\rho(\beta)U$; and if $\rho(\beta)U=U$, then $\hat\beta$
induces an automorphism of $E/U$ which acts on $(E/U)/\Phi(E/U)=
H/\Phi(H)$ as $\beta$ does.

\subsection{The Lazard correspondence}\label{subsec:lazard}
Let $p$ be odd. The Lazard correspondence \cite{Laz54},
\cite[Chapter~10]{Khu98} is an equivalence between finite $p$-groups of class
$<p$ and nilpotent Lie rings of $p$-power order and class $<p$; the two
structures live on the same set, the group operations being given by the
Baker--Campbell--Hausdorff formula and the Lie operations by its inverse
formulas, which are group words with exponents in $\Z_{(p)}$. A group of
exponent $p$ corresponds to a Lie algebra $\fl(G)$ over $K$ of the same
class. Automorphisms, normal subgroups and ideals, quotients, the lower
central series, and central subgroups and ideals correspond. We shall use
the correspondence only through Lemma~\ref{lem:lazard} below.

\subsection{Jennings determinacy}\label{subsec:det}

\begin{definition}\label{def:Ndet}
Let $N\ge2$ and let $G$ be a finite $p$-group with $D_N(G)=1$. We say
that $G$ is \emph{$N$-determined} if for every finite $p$-group $\Gt$ with
$K\Gt/\I(K\Gt)^N\cong KG/\I(KG)^N$ one has $\Gt/D_N(\Gt)\cong G$.
\end{definition}

If $G$ is $N$-determined it is $N'$-determined for all $N'\ge N$, since
$KG/\I^N$ is a quotient of $KG/\I^{N'}$; and then $G$ is determined by
$KG$, because $KG\cong K\Gt$ forces $\nu(\Gt)=\nu(G)<N$ by
Lemma~\ref{lem:grrigid}, so $D_N(\Gt)=1$. This is the form of the
conclusion in \cite[Theorem~3.1.2]{Roh90}, and it is the form that iterates:
the descent of Section~\ref{sec:descent} produces, from an isomorphism of
truncated algebras of two groups, an isomorphism of truncated algebras of
their bases, and it is truncated determinacy of the base, not
determinacy by the whole group algebra, that closes the induction.

\section{The free frame and the descent theorem}\label{sec:descent}

Throughout this section $p$ is an arbitrary prime.

\subsection{The truncated free frame}
Let $F$ be free on $x_1,\dots,x_d$. For $N\ge1$ let
\[
A_N:=K\langle X_1,\dots,X_d\rangle/(X_1,\dots,X_d)^N
\]
be the free associative algebra truncated in degree $N$, and let $\I_N$ be
its augmentation ideal, the image of $(X_1,\dots,X_d)$; it is a finite-dimensional local
algebra with $\I_N^{\,N}=0$ and $\gr A_N=K\langle X\rangle/(X)^N$. The
elements $1+X_i$ are units, and $x_i\mapsto1+X_i$ defines a surjective
algebra homomorphism $KF\to A_N$ carrying $\I(KF)^m$ onto $\I_N^{\,m}$.
For $R\lhd F$ of $p$-power index let $\fa_N(R)\subseteq\I_N$ be the ideal
generated by the images of the $r-1$, $r\in R$. Any assignment
$X_i\mapsto m_i\in\I_N$ extends uniquely to an algebra endomorphism of
$A_N$. The lifting principle in this truncated form is essentially
\cite[Lemma~3.1]{Roh87}, \cite[Proposition~3.1.1]{Roh90}.

\begin{lemma}\label{lem:frame}
Let $R\lhd F$ with $G=F/R$ a finite $p$-group, and $N\ge1$.
\begin{enumerate}[label=\textup{(\arabic*)}]
\item $A_N/\fa_N(R)\cong KG/\I(KG)^N$, compatibly with the maps from $KF$.
\item For $f\in F$, the image of $f-1$ lies in $\fa_N(R)$ if and only if
$f\in R\,M_N(F)$. In particular $\fa_N(R)=\fa_N(RM_N(F))$, and
$\fa_N(M_N(F))=0$, so that $A_N\cong K[F/M_N(F)]/\I^N$.
\item If $R\le\Phi(F)$ then $\fa_N(R)\subseteq\I_N^{\,2}$, and
$\I_N/\I_N^{\,2}\cong F/\Phi(F)\cong G/\Phi(G)$ when moreover $d(G)=d$.
\item An endomorphism $\Gamma$ of $A_N$ is an automorphism if and only if
the induced map $\gr_1\Gamma$ on $\I_N/\I_N^{\,2}$ is invertible; it is
then the identity on $\gr A_N$ if and only if $\gr_1\Gamma=\mathrm{id}$.
\end{enumerate}
\end{lemma}

\begin{proof}
(1) The map $KF\to A_N/\fa_N(R)$ kills $r-1$, hence factors through $KG$,
and kills $\I(KG)^N$; conversely $X_i\mapsto\bar x_i-1$ defines
$A_N\to KG/\I^N$ killing $\fa_N(R)$; the two maps are inverse. (2) By (1)
and Jennings' theorem, the image of $f-1$ lies in $\fa_N(R)$ iff
$\bar f-1\in\I(KG)^N$ iff $\bar f\in D_N(G)=M_N(F)R/R$. Since $u-1\in
\I(KF)^N$ for $u\in M_N(F)$, the image of $u-1$ is $0$. (3) For
$g,h\in F$, $g^p-1=(g-1)^p$ and $[g,h]-1=g^{-1}h^{-1}\bigl((g-1)(h-1)-
(h-1)(g-1)\bigr)$ lie in $\I(KF)^2$, and $xy-1=(x-1)+(y-1)+(x-1)(y-1)$;
so $s-1\in\I(KF)^2$ for $s\in\Phi(F)$. (4) $\gr A_N$ is generated in
degree one, so $\gr\Gamma$ is surjective in each degree if $\gr_1\Gamma$
is invertible, hence bijective, and then $\Gamma$ is injective on each
$\I_N^{\,m}/\I_N^{\,m+1}$, hence bijective; the converse is clear.
\end{proof}

We call $\Gamma\in\Aut(A_N)$ \emph{unipotent} if $\gr_1\Gamma=\mathrm{id}$,
that is, if $\Gamma(X_i)\equiv X_i$ modulo $\I_N^{\,2}$ for all $i$.

\begin{lemma}[normalized lift]\label{lem:lift}
Let $N\ge2$, let $G=F/R$ with $M_N(F)\subseteq R\subseteq\Phi(F)$ and
$d(G)=d$, let $\Gt$ be a finite $p$-group, and let
$\Theta:KG/\I^N\to K\Gt/\I^N$ be an algebra isomorphism. Then
$d(\Gt)=d$, and there are a presentation $\Gt=F/S$ with $S\le\Phi(F)$ and a
unipotent $\Gamma\in\Aut(A_N)$ with $\pi_S\circ\Gamma=\Theta\circ\pi_R$,
where $\pi_R:A_N\to KG/\I^N$ and $\pi_S:A_N\to K\Gt/\I^N$ are the maps of
Lemma~\ref{lem:frame}(1); consequently $\Gamma(\fa_N(R))=\fa_N(S)$.
Moreover, replacing $\Gt$ by $\Gt/D_N(\Gt)$ changes neither $K\Gt/\I^N$
nor $\Gt/D_N(\Gt)$, so one may assume $M_N(F)\subseteq S$.
\end{lemma}

\begin{proof}
$\Theta$ preserves the maximal ideals and induces
$\I/\I^2\cong\widetilde\I/\widetilde\I^{\,2}$, so $d(\Gt)=d$. Choose
$\tilde x_i\in\Gt$ with $\tilde x_i-1\equiv\Theta(x_i-1)$ modulo
$\widetilde\I^{\,2}$, possible because $g\mapsto g-1$ induces
$\Gt/\Phi(\Gt)\cong\widetilde\I/\widetilde\I^{\,2}$; by Burnside's basis
theorem the $\tilde x_i$ generate $\Gt$, and the kernel $S$ of
$x_i\mapsto\tilde x_i$ lies in $\Phi(F)$ because $d(\Gt)=d$. Define
$\Gamma$ by sending $X_i$ to a preimage in $\I_N$ of $\Theta(\pi_R(X_i))$
under $\pi_S$. Then $\pi_S\Gamma=\Theta\pi_R$ on generators, hence
everywhere, and $\Gamma(X_i)-X_i\in\pi_S^{-1}(\widetilde\I^{\,2})=
\fa_N(S)+\I_N^{\,2}=\I_N^{\,2}$ by Lemma~\ref{lem:frame}(3); so $\Gamma$
is unipotent and an automorphism. Since $\Theta$ is injective,
$\Gamma(\fa_N(R))=\Gamma(\ker\pi_R)\subseteq\ker\pi_S=\fa_N(S)$, and the two ideals have the same
codimension, so they are equal. The last sentence holds because
$a-1\in\I(K\Gt)^N$ for $a\in D_N(\Gt)$, so $K\Gt/\I^N=K[\Gt/D_N(\Gt)]/\I^N$.
\end{proof}

\begin{lemma}[free transfer of the top Jennings layer]\label{lem:transfer}
Let $n\ge1$. Every automorphism $\Gamma$ of $A_{n+1}$ satisfies
$\Gamma(\fa_{n+1}(M_n(F)))=\fa_{n+1}(M_n(F))$. Consequently, if
$R,S\lhd F$ have $p$-power index and $\Gamma(\fa_{n+1}(R))=\fa_{n+1}(S)$,
then $\Gamma(\fa_{n+1}(RM_n(F)))=\fa_{n+1}(SM_n(F))$.
\end{lemma}

\begin{proof}
By Lemma~\ref{lem:frame}(2), $A_{n+1}=KP/\I(KP)^{n+1}$ for the finite
$p$-group $P=F/M_{n+1}(F)$, so $\gr A_{n+1}=u(L(P))$ in degrees $\le n$
by Quillen's theorem. As $\Gamma$ preserves the radical $\I_{n+1}$ and
its powers, $\gr\Gamma$ is a graded automorphism; it preserves
$L_1(P)=\I/\I^2$, hence the brackets and $p$-th powers of degree-one
classes that land in degrees $\le n$, hence each $L_m(P)$ for $m\le n$,
since $L(P)$ is generated in degree one. Now $L_n(P)$ is the image of
$M_n(P)=M_n(F)/M_{n+1}(F)$, so for $u\in M_n(F)$ the image of $u-1$ in
$\I_{n+1}^{\,n}=\I_{n+1}^{\,n}/\I_{n+1}^{\,n+1}$ lies in $L_n(P)$, and
so does its image under $\Gamma$: there is $u'\in M_n(F)$ with
$\Gamma(u-1)=u'-1$ in $A_{n+1}$. Hence $\Gamma(\fa_{n+1}(M_n(F)))
\subseteq\fa_{n+1}(M_n(F))$, and equality holds by dimension. The second
statement follows from $\fa_{n+1}(RM_n(F))=\fa_{n+1}(R)+
\fa_{n+1}(M_n(F))$.
\end{proof}

\begin{remark}\label{rem:intrinsic}
Downstairs the lemma says: if $A=D_n(G)$ then $\I(KA)KG+\I^{\,n+1}=
\operatorname{span}\{a-1:a\in A\}+\I^{\,n+1}$, because $(a-1)\I\subseteq
\I^{\,n+1}$; and this is the preimage of $L_n(G)\subseteq\I^n/\I^{n+1}$,
an invariant of the algebra by Lemma~\ref{lem:grrigid}. So the ideal
generated by the top Jennings layer transfers along every algebra
isomorphism \emph{modulo $\I^{\,n+1}$}, with no hypothesis. It need not
transfer on the nose, and by Sandling's example \cite[\S3.2]{Roh90} no
general correspondence of ideals is available; Lemma~\ref{lem:U} below
says that, for the purpose of recovering $G$, the transfer modulo
$\I^{\,n+1}$ suffices.
\end{remark}

\subsection{The action on the relation module}\label{subsec:action}
Fix $n\ge2$ and $R\lhd F$ with
\begin{equation}\label{eq:standing}
M_{n+1}(F)\subseteq R\subseteq\Phi(F),\qquad d(F/R)=d ,
\end{equation}
and put $G=F/R$ and $A=D_n(G)$, assumed non-trivial, so that
$n=\nu(G)$. Put $R'=R\,M_n(F)$ and $H=F/R'=G/A$. Since $[A,G]\le
D_{n+1}(G)=1$ and $A^p\le D_{np}(G)=1$, $A$ is central and elementary
abelian, and $A\le D_2(G)=\Phi(G)$; so $A\to G\to H$ is a central
Frattini extension and $d(H)=d$. Let $L'=[F,R'](R')^p$, $E=F/L'$,
$V=R'/L'$, $W_k=V\cap D_k(E)$ be the Hertweck--Soriano data of $H$ built
on the presentation $F/R'$, and put $U=R/L'$.

\begin{lemma}\label{lem:U}
$L'\subseteq R$, so $U\le V$ and $G=E/U$ with $A=V/U$; and
$W_{n+1}\subseteq U$.
\end{lemma}

\begin{proof}
$[F,R']R/R=[G,A]=1$ and $(R')^pR/R=A^p=1$, so $L'\subseteq R$; and
$V/U=R'/R=A$. The image of $D_{n+1}(E)$ in $G=E/U$ is $D_{n+1}(G)=1$, so
$D_{n+1}(E)\subseteq U$.
\end{proof}

Write $\fa=\fa_{n+1}$ and $\I=\I_{n+1}$ in $A_{n+1}$, and put
\[
\fI:=\fa(R')\,\I+\I\,\fa(R'),\qquad \st:=A_{n+1}/\fI .
\]

\begin{lemma}\label{lem:st}
$\st\cong s(E,V)/\I(s)^{n+1}$, compatibly with the maps from $KF$.
\end{lemma}

\begin{proof}
The kernel of $KF\to s(E,V)$ is $\Delta(F,L')+\Delta(F,R')\I(KF)$, where
$\Delta(F,N)$ denotes the ideal generated by the $r-1$, $r\in N$. For
$r\in R'$ and $f\in F$ one has $r^p-1=(r-1)^p\in\Delta(F,R')\I(KF)$ and
\[
[r,f]-1=r^{-1}f^{-1}\bigl((r-1)(f-1)-(f-1)(r-1)\bigr),
\]
so $\Delta(F,L')\subseteq\Delta(F,R')\I(KF)+\I(KF)\Delta(F,R')$;
conversely the same identity gives $(f-1)(r-1)\in\Delta(F,R')\I(KF)+
\Delta(F,L')$. Hence the kernel is $\Delta(F,R')\I(KF)+\I(KF)\Delta(F,R')$
(cf.\ \cite[Lemma~2.6]{Roh87}), and the kernel of $KF\to s(E,V)/\I(s)^{n+1}$
is obtained by adding $\I(KF)^{n+1}$; its image in $A_{n+1}$ is $\fI$.
\end{proof}

\begin{definition}\label{def:Autn}
$\Aut_n(H):=\{\Gamma\in\Aut(A_{n+1}):\Gamma(\fa(R'))=\fa(R')\}$, and
$\UAut_n(H)$ is its subgroup of unipotent elements. By
Lemma~\ref{lem:frame}(1),(4) and the universal property of $A_{n+1}$,
every automorphism of $KH/\I^{\,n+1}$ is induced by some element of
$\Aut_n(H)$.
\end{definition}

\begin{lemma}\label{lem:rho}
\begin{enumerate}[label=\textup{(\arabic*)}]
\item Every $\Gamma\in\Aut_n(H)$ preserves $\fI$, hence induces
$\tilde\Gamma\in\Aut(\st)$, and $\tilde\Gamma$ preserves the image
$\widetilde{\mathcal V}$ of $\fa(R')$ in $\st$. If $\Gamma$ is unipotent
then $(\tilde\Gamma-1)\I(\st)^m\subseteq\I(\st)^{m+1}$ for all $m$.
\item Under Lemma~\ref{lem:st}, $\widetilde{\mathcal V}$ is the image of
$\{V-1\}$, and $v\mapsto v-1$ induces an isomorphism of $K$-spaces
$\Vt:=V/W_{n+1}\to\widetilde{\mathcal V}$. Transporting
$\tilde\Gamma|_{\widetilde{\mathcal V}}$ gives a homomorphism
$\rho_n:\Aut_n(H)\to GL(\Vt)$, and $\rho_n(\Gamma)$ preserves the images
of the $W_k$.
\item Let $\beta\in\Aut(H)$, lifted to $\theta\in\Aut(F)$ with
$\theta(R')=R'$ as in \S\ref{subsec:HS}, and let $\hat\theta$ be the
automorphism of $A_{n+1}$ defined by $X_i\mapsto\theta(x_i)-1$. Then
$\hat\theta\in\Aut_n(H)$ and $\rho_n(\hat\theta)$ is the map induced by
$\rho(\beta)$ on $V/W_{n+1}$.
\item For $S\lhd F$ with $L'\subseteq S\subseteq R'$, the image of
$\fa(S)$ in $\st$ corresponds to $(S/L')W_{n+1}/W_{n+1}\subseteq\Vt$. In
particular the image of $\fa(R)$ corresponds to $\widetilde U:=U/W_{n+1}$.
\end{enumerate}
\end{lemma}

\begin{proof}
(1) $\Gamma(\fa(R')\I)=\fa(R')\I$ since $\Gamma(\I)=\I$, and likewise on
the other side; the same for $\Gamma^{-1}$. For unipotent $\Gamma$,
$\tilde\Gamma(X_i)-X_i\in\I^2$, and the telescoping identity
\[
\tilde\Gamma(a_1\cdots a_m)-a_1\cdots a_m=\sum_j\tilde\Gamma(a_1)\cdots
\tilde\Gamma(a_{j-1})\bigl(\tilde\Gamma(a_j)-a_j\bigr)a_{j+1}\cdots a_m
\]
gives the depth statement. (2) The image of $\fa(R')$ in $s(E,V)$ is the
image of $\I(KV)KE$, that is $\{V-1\}$, and $\{V-1\}\cap\I(s)^{n+1}=
\{W_{n+1}-1\}$ by Lemma~\ref{lem:PS}; additivity of $v\mapsto v-1$ is
\cite[(3.1)]{HS07}. (3) $\hat\theta$ maps $r-1$ to $\theta(r)-1$, so it
preserves $\fa(R')$, and on $s(E,V)$ it is induced by the automorphism of
$E$ lifting $\beta$. (4) In $s(E,V)$ the ideal generated by
$\{s-1:s\in S/L'\}$ is their span, by \eqref{eq:trunc0} and additivity,
namely $\{S/L'-1\}$.
\end{proof}

\subsection{The descent theorem}

\begin{theorem}[descent]\label{thm:descent}
Let $p$ be any prime, $n\ge2$, and $R\lhd F$ as in \eqref{eq:standing},
with $G=F/R$, $A=D_n(G)\ne1$, $H=G/A$, and $V,U,\Vt,\widetilde U$ as in
\S\ref{subsec:action}. Assume:
\begin{enumerate}[label=\textup{(\alph*)}]
\item $H$ is $(n+1)$-determined;
\item $\rho_n(\Aut_n(H))\subseteq\rho(\Aut(H))|_{\Vt}$.
\end{enumerate}
Then $G$ is $(n+1)$-determined.
\end{theorem}

\begin{proof}
Let $\Gt$ be a finite $p$-group with $K\Gt/\I^{\,n+1}\cong KG/\I^{\,n+1}$.
By Lemma~\ref{lem:lift} we may assume $D_{n+1}(\Gt)=1$ and fix
$\Gt=F/S$ with $M_{n+1}(F)\subseteq S\subseteq\Phi(F)$ and a unipotent
$\Gamma\in\Aut(A_{n+1})$ with $\Gamma(\fa(R))=\fa(S)$.

\emph{Step 1: the bases agree.} Put $S'=SM_n(F)$ and $\widetilde H=F/S'$.
By Lemma~\ref{lem:transfer}, $\Gamma(\fa(R'))=\fa(S')$, which by
Lemma~\ref{lem:frame}(1) is an isomorphism $KH/\I^{\,n+1}\cong
K\widetilde H/\I^{\,n+1}$; and $D_{n+1}(\widetilde H)=1$. Hypothesis (a)
gives $\widetilde H\cong H$.

\emph{Step 2: normalization.} An isomorphism $F/R'\to F/S'$ lifts to
$\theta\in\operatorname{End}(F)$ with $\theta(R')\subseteq S'$; $\theta$ is
surjective modulo $\Phi(F)$, hence surjective, hence an automorphism, and
$\theta(R')=S'$ by comparing orders. Put $S_1=\theta^{-1}(S)$ and
$\Gamma_1=\hat\theta^{-1}\Gamma$. Then $\Gamma_1(\fa(R))=\fa(S_1)$,
$S_1M_n(F)=\theta^{-1}(S')=R'$ as $M_n(F)$ is characteristic, so
$\Gamma_1\in\Aut_n(H)$ by Lemma~\ref{lem:transfer}; also $F/S_1\cong\Gt$,
and $S_1$ satisfies \eqref{eq:standing} with $D_n(F/S_1)=R'/S_1$, so
Lemma~\ref{lem:U} applied to $S_1$ gives $L'\subseteq S_1$ and
$W_{n+1}\subseteq U_{S_1}:=S_1/L'$.

\emph{Step 3: the subspaces.} By Lemma~\ref{lem:rho}(1),(4),
$\tilde\Gamma_1$ carries the image of $\fa(R)$, which corresponds to
$\widetilde U$, onto the image of $\fa(S_1)$, which corresponds to
$\widetilde U_{S_1}=U_{S_1}/W_{n+1}$: thus
$\rho_n(\Gamma_1)\widetilde U=\widetilde U_{S_1}$.

\emph{Step 4: the group.} By (b) there is $\beta\in\Aut(H)$ with
$\rho_n(\Gamma_1)=\rho(\beta)|_{\Vt}$. As $\rho(\beta)$ preserves
$W_{n+1}$ and both $U$ and $U_{S_1}$ contain $W_{n+1}$, this lifts to
$\rho(\beta)U=U_{S_1}$, and the easy half of the Hertweck--Soriano
correspondence (\S\ref{subsec:HS}) gives $\Gt\cong E/U_{S_1}\cong E/U=G$.
\end{proof}

When hypothesis (b) fails, the proof still bounds the failure.

\begin{proposition}[counting]\label{prop:count}
Let $p$ be any prime, $n\ge2$, and $R\lhd F$ as in \eqref{eq:standing},
with $G=F/R$, $A=D_n(G)\ne1$, $H=G/A$, and $\Vt$ as in
\S\ref{subsec:action}. Assume that $H$ is $(n+1)$-determined, and put
\[
\mathcal G=\rho_n(\Aut_n(H)),\qquad
\mathcal A=\rho(\Aut(H))|_{\Vt}\qquad(\le GL(\Vt)).
\]
Then $\mathcal A\le\mathcal G$, and the number of isomorphism classes of
finite $p$-groups $\Gt$ with $D_{n+1}(\Gt)=1$ and
$\Fp\Gt/\I^{\,n+1}\cong\Fp G/\I^{\,n+1}$ is at most
$[\mathcal G:\mathcal A]$. Hypothesis \textup{(b)} of
Theorem~\ref{thm:descent} is the case $\mathcal G=\mathcal A$.
\end{proposition}

\begin{proof}
By Lemma~\ref{lem:rho}(3), for $\beta\in\Aut(H)$ the automorphism
$\hat\theta$ lies in $\Aut_n(H)$ and $\rho_n(\hat\theta)=\rho(\beta)|_{\Vt}$;
so $\mathcal A\le\mathcal G$, both being groups because $\rho_n$ and $\rho$
are homomorphisms. For $\gamma\in\mathcal G$ let $U_\gamma\le V$ be the
preimage of $\gamma\widetilde U$; it contains $W_{n+1}$. Steps~1--3 of the
proof of Theorem~\ref{thm:descent}, which use only hypothesis (a), show
that every $\Gt$ as in the statement is isomorphic to
$E/U_{S_1}=E/U_\gamma$ with $\gamma=\rho_n(\Gamma_1)\in\mathcal G$. If
$\gamma'=\rho(\beta)|_{\Vt}\,\gamma$ with $\beta\in\Aut(H)$, then
$U_{\gamma'}=\rho(\beta)U_\gamma$, because $\rho(\beta)$ preserves $W_{n+1}$;
so $E/U_{\gamma'}\cong E/U_\gamma$ by the easy half of the
Hertweck--Soriano correspondence (\S\ref{subsec:HS}). Hence the isomorphism
class of $E/U_\gamma$ depends only on the coset $\mathcal A\gamma$, and every
class in question arises in this way.
\end{proof}

Hypothesis (b) is the one that carries the mathematics. We split it.

\begin{definition}\label{def:props}
Let $X=F/R_X$ be a finite $p$-group with $d(X)=d$.
\begin{enumerate}[label=\textup{(\roman*)}]
\item For $m\ge1$ put $\Aut^{(m)}(R_X):=\{\Gamma\in\Aut(A_m):
\Gamma(\fa_m(R_X))=\fa_m(R_X)\}$ and let $\overline{\Aut}^{(m)}(X)$ be its
image in $GL(\I_m/\I_m^{\,2})=GL(X/\Phi(X))$. We say that $X$ has property
$(\Lev_m)$ if $\overline{\Aut}^{(m+1)}(X)$ is contained in the image of
$\Aut(X)$ in $GL(X/\Phi(X))$.
\item For $n\ge2$, with $H=X$ in the notation of \S\ref{subsec:action},
we say that $H$ has property $(\RR_n)$ (\emph{truncated relative rigidity
at level $n$}) if for every $\Gamma\in\UAut_n(H)$ there is
$\beta\in\Aut(H)$, acting trivially on $H/\Phi(H)$, with
$\rho_n(\Gamma)=\rho(\beta)|_{\Vt}$.
\end{enumerate}
\end{definition}

Note that $\Aut_n(H)=\Aut^{(n+1)}(R')$ and that the notions depend on $X$
only, since $A_m$ and $\fa_m$ depend on the presentation only through the
choice of the free generators, which conjugates everything.

\begin{lemma}\label{lem:split}
Hypothesis \textup{(b)} of Theorem~\ref{thm:descent} holds if $H$ has
$(\Lev_{n-1})$ and $(\RR_n)$.
\end{lemma}

\begin{proof}
Let $\Gamma\in\Aut_n(H)$, so that $\Gamma\in\Aut^{(n+1)}(R')$. Reduction
modulo $\I^{\,n}$
maps $\Aut^{(n+1)}(R')$ into $\Aut^{(n)}(R')$ without changing $\gr_1$, so
$(\Lev_{n-1})$ provides $\beta_0\in\Aut(H)$ with the same image in
$GL(H/\Phi(H))$ as $\Gamma$. Lift $\beta_0$ to $\theta_0\in\Aut(F)$ with
$\theta_0(R')=R'$; then $\hat\theta_0^{-1}\Gamma\in\Aut_n(H)$ is
unipotent, and $(\RR_n)$ gives $\beta_1$ with
$\rho_n(\hat\theta_0^{-1}\Gamma)=\rho(\beta_1)|_{\Vt}$. Hence
$\rho_n(\Gamma)=\rho(\beta_0\beta_1)|_{\Vt}$ by Lemma~\ref{lem:rho}(3).
\end{proof}

\begin{lemma}\label{lem:monotone}
Let $X=F/R_X$ with $d(X)=d$.
\begin{enumerate}[label=\textup{(\arabic*)}]
\item $(\Lev_m)$ implies $(\Lev_{m'})$ for all $m'\ge m$.
\item If $X$ has $(\Lev_m)$ for some $m$, then
$\Aut(KX)=\Aut(X)\cdot\UAut(KX)$.
\item If $X$ is relatively free, that is $R_X$ is fully invariant in $F$,
then $X$ has $(\Lev_m)$ for all $m$.
\end{enumerate}
\end{lemma}

\begin{proof}
(1) Reduction modulo $\I^{\,m+1}$ maps $\Aut^{(m'+1)}(R_X)$ into
$\Aut^{(m+1)}(R_X)$, preserving $\gr_1$. (2) An automorphism $\theta$ of
$KX$ induces one of $KX/\I^{\,m+1}$, which lifts as in the proof of
Lemma~\ref{lem:lift} (without normalization) to $\Gamma\in\Aut(A_{m+1})$
with $\Gamma(\fa_{m+1}(R_X))\subseteq\fa_{m+1}(R_X)$, hence equal, and
with $\gr_1\Gamma=\gr_1\theta$. By $(\Lev_m)$ there is $\beta\in\Aut(X)$
inducing $\gr_1\theta$, and $\beta^{-1}\theta\in\UAut(KX)$. (3) If $X$ is
relatively free, every linear automorphism of $X/\Phi(X)$ lifts to an
endomorphism of $X$, which is surjective by Burnside's theorem; so the
image of $\Aut(X)$ is all of $GL(X/\Phi(X))$.
\end{proof}

\begin{theorem}[propagation of the Levi factorization]\label{thm:levi}
Let $p$ be any prime, $n\ge2$, $R\lhd F$ as in \eqref{eq:standing},
$G=F/R$ with $A=D_n(G)\ne1$ and $H=G/A$. If $H$ has $(\Lev_{n-1})$ and
$(\RR_n)$, then $G$ has $(\Lev_n)$.
\end{theorem}

\begin{proof}
Let $\Gamma\in\Aut^{(n+1)}(R)$ with image $g\in GL(G/\Phi(G))=
GL(H/\Phi(H))$. By Lemma~\ref{lem:transfer}, $\Gamma\in\Aut_n(H)$, and as
in the proof of Lemma~\ref{lem:split} we find $\beta_0,\beta_1\in\Aut(H)$,
with $\beta_1$ trivial on $H/\Phi(H)$ and $\beta_0$ inducing $g$, such
that $\rho_n(\Gamma)=\rho(\beta_0\beta_1)|_{\Vt}$. On the other hand
$\Gamma$ preserves $\fa(R)$, so $\rho_n(\Gamma)\widetilde U=\widetilde U$
by Lemma~\ref{lem:rho}(4), whence $\rho(\beta_0\beta_1)U=U$ as
$W_{n+1}\subseteq U$. By \S\ref{subsec:HS}, $\beta_0\beta_1$ lifts to an
automorphism of $G=E/U$ inducing $\beta_0\beta_1$, that is $g$, on
$G/\Phi(G)$.
\end{proof}

\begin{corollary}[R\"ohl's theorem]\label{cor:rohl}
Let $p$ be any prime and $R\lhd F$ with $M_{n+1}(F)\subseteq R\subsetneq
M_n(F)$. Then $F/R$ is $(n+1)$-determined; in particular it is determined
by its modular group algebra \cite[Theorem~3.1.2]{Roh90}.
\end{corollary}

\begin{proof}
Induction on $n$, the case $n=1$ being that of an elementary abelian
group, which is $2$-determined since $KX/\I^2$ determines $d(X)$. For
$G=F/R$ the base is $H=F/M_n(F)$, relatively free, so it has
$(\Lev_{n-1})$ by Lemma~\ref{lem:monotone}(3) and is $n$-determined, hence
$(n+1)$-determined, by induction (its relations $M_n(F)$ lie between
$M_n$ and $M_{n-1}$). As $R'=M_n(F)$, one has $V=W_n$, and for unipotent
$\Gamma\in\UAut_n(H)$ the depth statement of Lemma~\ref{lem:rho}(1) gives
$\tilde\Gamma(v-1)-(v-1)\in\I(\st)^{n+1}=0$ for $v\in V$: so
$\rho_n(\Gamma)=1$ and $(\RR_n)$ holds with $\beta=1$. Apply
Lemma~\ref{lem:split} and Theorem~\ref{thm:descent}.
\end{proof}

So R\"ohl's hypothesis is the case in which the relation module sits in a
single layer of the depth flag, where a unipotent automorphism cannot
move it. When the relations occupy two or more layers, $(\RR_n)$ becomes a
genuine condition, and it is the condition that Section~\ref{sec:rigid}
establishes. Another proof of R\"ohl's theorem through small group rings
is given in \cite[Theorem~5.8]{HS07}.

\section{Truncated relative rigidity in exponent $p$}\label{sec:rigid}

Throughout this section $p$ is odd and $H$ is a finite $p$-group of
exponent $p$ and class $c\le p-2$; we keep the notation of
\S\ref{subsec:HS} and \S\ref{subsec:action} for $H$, with $R'\le\Phi(F)$
its relation subgroup, $E$, $V$, $W_k$. Since $H$ has exponent $p$,
$F^{\,p}\subseteq R'$. Put
\[
C(V):=V\cap[E,E],\qquad E^\circ:=E/\mho_1(E),\qquad
V^\circ:=V/\mho_1(E).
\]

The results of this section are truncated analogues of results in the
author's paper \cite{UA}, where the untruncated action of $\UAut(KH)$ on
$V$ is treated at truncation degree $p$; the shadow argument of
Lemma~\ref{lem:shadow} appears there. Since the
statements needed here are different --- levels $n\le p-1$ in
Theorem~\ref{thm:trunc}, and the reduction from the ceiling $n=p$ in
Theorem~\ref{thm:ceiling}, which the untruncated theory does not provide
--- and since the paper should be readable on its own, complete proofs are
given.

\subsection{The universal extension of a base of exponent $p$}

\begin{proposition}\label{prop:E}
\begin{enumerate}[label=\textup{(\arabic*)}]
\item $E$ has class at most $c+1\le p-1$, $\gamma_2(E)$ has exponent
$p$, and $g\mapsto g^p$ is a homomorphism $E\to V$ which factors through
$E/\Phi(E)\cong H/\Phi(H)$; its image $\mho_1(E)$ is elementary abelian of
rank $d$, spanned by the classes of $x_1^{\,p},\dots,x_d^{\,p}$.
\item $E/[E,E]\cong(C_{p^2})^d$, generated by the images of the $x_i$,
and the natural map $\varphi:E\to E/[E,E]$ is injective on $\mho_1(E)$
with $\varphi(V)=\varphi(\mho_1(E))$; hence
\[
V=C(V)\oplus\mho_1(E),\qquad V^\circ\cong C(V).
\]
\item $D_m(E)=\gamma_m(E)\mho_1(E)$ for $m\le p$ and $D_m(E)=1$ for
$m\ge p+1$. Consequently $W_{p+1}=1$, and $W_m=\mho_1(E)$ for
$c+1<m\le p$; in particular $C(V)\cap W_p=0$.
\item $E^\circ$ has exponent $p$ and class at most $c+1$, $V^\circ$ is
central in $E^\circ$ and contained in $\Phi(E^\circ)=\gamma_2(E^\circ)$,
and $E^\circ/V^\circ=H$.
\end{enumerate}
\end{proposition}

\begin{proof}
(1) $E/V=H$ has class $c$ and $V$ is central, so $\gamma_{c+2}(E)=1$;
and $g^p\in V$ for all $g\in E$. We use the Hall--Petrescu formula
\cite[III.9.4]{Hup67}: for $x,y$ in a group and $m\ge1$,
\begin{equation}\label{eq:HP}
x^my^m=(xy)^m\,c_2^{\binom m2}c_3^{\binom m3}\cdots
c_{m-1}^{\binom m{m-1}}c_m,\qquad c_i\in\gamma_i(\langle x,y\rangle).
\end{equation}
We show by downward induction on $k$ that $\gamma_k(E)$ has exponent $p$
for all $k\ge2$; this is clear for $k\ge p$, as $\gamma_p(E)=1$. Assume it
for $k+1$, where $2\le k\le p-1$. For $g,h\in\gamma_k(E)$, the $c_i$ in
\eqref{eq:HP} with $m=p$ lie in $\gamma_{ik}(E)\subseteq\gamma_{k+1}(E)$
for $i\ge2$, so $c_i^{\binom pi}=1$ for $2\le i\le p-1$ and $c_p\in
\gamma_p(E)=1$: thus $(gh)^p=g^ph^p$ on $\gamma_k(E)$. Since $\gamma_k(E)$
is generated by the $[a,b]$ with $a\in\gamma_{k-1}(E)$, $b\in E$, it
suffices to show $[a,b]^p=1$. Apply \eqref{eq:HP} to $x=a$, $y=[a,b]$: a
simple commutator of weight $i\ge2$ in $a$ and $[a,b]$ is trivial unless
both occur, in which case it lies in $\gamma_{2k-1}(E)\subseteq
\gamma_{k+1}(E)$; hence again all correction terms are trivial and
$a^p[a,b]^p=(a[a,b])^p=(a^b)^p=(a^p)^b=a^p$, the last step because
$a^p\in V$ is central. So $[a,b]^p=1$. Now for arbitrary $g,h\in E$ the
$c_i$ lie in $\gamma_2(E)$, which has exponent $p$, so $(gh)^p=g^ph^p$:
the $p$-th power map is a homomorphism $E\to V$. It kills $\gamma_2(E)$
and $E^p\subseteq V$, hence factors through $E/\Phi(E)$, and its image is
spanned by the $x_i^{\,p}$, whose independence follows from (2).

(2) $E/[E,E]=F/[F,F]L'=F/[F,F](R')^p$. As $F^{\,p}\subseteq R'\subseteq
F^{\,p}[F,F]$, the image of $R'$ in $F/[F,F]\cong\Z^d$ is $p\Z^d$ and that
of $(R')^p$ is $p^2\Z^d$; so $E/[E,E]\cong(C_{p^2})^d$, the images of the
$x_i^{\,p}$ are independent of order $p$, and $\varphi(V)=p\Z^d/p^2\Z^d$
is spanned by them. Thus $\varphi$ is injective on $\mho_1(E)$ with
$\varphi(\mho_1(E))=\varphi(V)$, and $V=C(V)\oplus\mho_1(E)$ with
$C(V)=\ker(\varphi|_V)$.

(3) In \eqref{eq:jennings} for $E$, the factors $\mho_j(\gamma_i(E))$
with $i\ge2$, $j\ge1$ are trivial by (1), and $\mho_2(E)=1$ since
$\mho_1(E)\le V$ is elementary abelian; so $D_m(E)=\gamma_m(E)\mho_1(E)$
for $m\le p$ and $D_m(E)=\gamma_m(E)=1$ for $m\ge p+1$. The statements
on $W_m$ follow from $\mho_1(E)\le V$, the modular law, and
$\gamma_m(E)=1$ for $m\ge c+2$; and $C(V)\cap W_p=C(V)\cap\mho_1(E)=0$.

(4) $E^\circ$ is a quotient of $E$, has exponent $p$ since $g^p\in
\mho_1(E)$, and $V^\circ$ is central with $E^\circ/V^\circ=E/V=H$; also
$V\le\Phi(E)$, and $\Phi(E^\circ)=\gamma_2(E^\circ)$ as $E^\circ$ has
exponent $p$.
\end{proof}

\begin{lemma}\label{lem:liftsE}
Let $\beta\in\Aut(H)$ and $\hat\beta\in\Aut(E)$ a lift as in
\S\ref{subsec:HS}. Then $\hat\beta$ preserves $[E,E]$, $C(V)$ and
$\mho_1(E)$, hence descends to an automorphism $\beta^\circ$ of $E^\circ$
preserving $V^\circ$ and inducing $\beta$ on $H$; and $\rho(\beta)$ acts on
$V^\circ\cong C(V)$ as $\beta^\circ$ does. If $\beta$ acts trivially on
$H/\Phi(H)$, then $\hat\beta$ fixes $\mho_1(E)$ pointwise.
\end{lemma}

\begin{proof}
The first assertions are clear. If $\beta$ is trivial on $H/\Phi(H)=
E/\Phi(E)$, then $\hat\beta(x)\in x\Phi(E)$ and, the $p$-th power map
being a homomorphism killing $\Phi(E)$, $\hat\beta(x^p)=\hat\beta(x)^p=
x^p$.
\end{proof}

\subsection{Logarithms}

\begin{lemma}[logarithm of a unipotent automorphism]\label{lem:log}
Let $A$ be a finite-dimensional $K$-algebra with an ideal $I$ such that
$A=K\cdot1\oplus I$ and $I^{\,p}=0$, and let $\alpha\in\Aut(A)$ satisfy
$(\alpha-1)(I)\subseteq I^2$. Put $N=\alpha-1$ and
$\Delta=\sum_{k=1}^{p-2}\frac{(-1)^{k+1}}{k}N^k$. Then $\Delta$ is a
derivation of $A$ with $\Delta(I^m)\subseteq I^{m+1}$, $\Delta^{p-1}=0$,
and $\alpha=\exp\Delta:=\sum_{k=0}^{p-2}\Delta^k/k!$.
\end{lemma}

\begin{proof}
Telescoping as in Lemma~\ref{lem:rho}(1) gives $N(I^m)\subseteq I^{m+1}$;
hence $N^{p-1}=0$ and $\Delta(I^m)\subseteq I^{m+1}$, $\Delta^{p-1}=0$.
In $\Z_{(p)}[t]/(t^{p-1})$ the truncated identities
$\exp(\log(1+t))=1+t$ and $\log(\exp t)=t$ hold, all denominators being
prime to $p$; substituting $t\mapsto N$ gives $\alpha=\exp\Delta$.

For the Leibniz rule let $Q$ be the quotient of $I\otimes I$ by
$\sum_{i+j\ge p}I^i\otimes I^j$. The multiplication $m:I\otimes I\to A$
factors through $Q$, and the commuting operators $X=\Delta\otimes1$,
$Y=1\otimes\Delta$ act on $Q$ raising total degree by one, so every
monomial $X^aY^b$ with $a+b\ge p-2$ vanishes on $Q$. Hence on $Q$
\[
\alpha\otimes\alpha=\exp(X)\exp(Y)=\exp(X+Y),
\]
the last equality because $\exp(X)\exp(Y)-\exp(X+Y)$ is a combination
of monomials of total degree $\ge p-1$ in $\Z_{(p)}[X,Y]$. Put
$M'=\alpha\otimes\alpha-1$ on $Q$; then $(X+Y)^{p-2}=0$ on $Q$, so
$\log(1+M')=X+Y$ there. Since $m\circ(\alpha\otimes\alpha)=\alpha\circ m$,
we get $m\circ M'^k=N^k\circ m$ for all $k$, hence
$m\circ\log(1+M')=\Delta\circ m$, that is
$\Delta(ab)=\Delta(a)b+a\Delta(b)$ for $a,b\in I$; and trivially when one
of them is a scalar.
\end{proof}

\begin{lemma}[the Lazard algebra inside the truncated group algebra]
\label{lem:lazard}
Let $p$ be odd, $G$ a finite $p$-group of exponent $p$ and class $<N$,
where $N\le p$, and $B=KG/\I^N$. Then $\exp(x)=\sum_{k<p}x^k/k!$ and
$\log(1+x)=\sum_{1\le k<p}(-1)^{k+1}x^k/k$ are mutually inverse bijections
$\I(B)\to1+\I(B)$ with $\exp(x)\exp(y)=\exp(Z(x,y))$, where $Z$ is the
Baker--Campbell--Hausdorff series truncated in degree $p$. The group $G$
embeds in $B^\times$, and $\fl:=\log G$ is a Lie subalgebra of
$(B,[\ ,\ ])$ with $\fl\cap\I(B)^m=\gamma_m(\fl)$ for all $m$. The
bijection $\log:G\to\fl$ realizes the Lazard correspondence: the Lie
algebra $\fl$ is $\fl(G)$, normal subgroups $M$ of $G$ correspond to the
ideals $\log M$, central subgroups to central ideals, quotients to
quotients, and a Lie automorphism $\Lambda$ of $\fl$ corresponds to the
group automorphism $\exp\circ\Lambda\circ\log$ of $G$.
\end{lemma}

\begin{proof}
The identities for $\exp$, $\log$ and $Z$ hold in the free associative
power series algebra over $\Z_{(p)}$ modulo degree $p$, with denominators
prime to $p$, hence in $B$. $D_N(G)=1$ by \eqref{eq:jennings}, so $G$
embeds in $B^\times$. The Lie ring $(\I(B),[\ ,\ ])$ is nilpotent of
class $<p$ and of $p$-power order, so by the Lazard correspondence
\cite[Chapter~10]{Khu98} the operation $Z$ makes it a group and the inverse
Baker--Campbell--Hausdorff formulas, which are group words with exponents
in $\Z_{(p)}$, recover its Lie operations. Since $\log(gh)=Z(\log g,\log
h)$ and $\log(g^k)=k\log g$, the set $\fl=\log G$ is a subgroup of
$(\I(B),Z)$, hence closed under the Lie operations of $\I(B)$: it is a Lie
subalgebra, and $(\fl,Z)\cong G$ via $\exp$. So $\fl$ is the Lazard Lie
algebra of $G$, and the correspondence of subgroups, ideals and
automorphisms is that of the general theory. Finally $\log g=(g-1)u$ with
$u\in1+\I(B)$ commuting with $g-1$, so $\log g\in\I(B)^m$ if and only if
$g\in D_m(G)=\gamma_m(G)$, and $\log\gamma_m(G)=\gamma_m(\fl)$.
\end{proof}

\subsection{The symmetrization decomposition}
Apply Lemma~\ref{lem:lazard} to $G=E^\circ$, which has exponent $p$ and
class at most $c+1$ (Proposition~\ref{prop:E}(4)), with any $N$ such that
$c+1<N\le p$: put
\[
B:=KE^\circ/\I^N,\qquad\fe:=\log E^\circ\subseteq B,\qquad
\fv:=\log V^\circ\subseteq\fe .
\]
Thus $\fe=\fl(E^\circ)$, $\fv$ is a central ideal of $\fe$ contained in
$\gamma_2(\fe)$ (as $V^\circ\le\gamma_2(E^\circ)$ is central), and
$\fe/\fv=\fl(H)$. For $2\le j\le N-1$ let $\varepsilon_j:S^j\fe\to B$ be the
symmetrization
\[
\varepsilon_j(u_1\cdots u_j)=\frac1{j!}\sum_{\pi\in S_j}
u_{\pi(1)}\cdots u_{\pi(j)}.
\]
For $b\in\fe$, $\ad_b=[b,-]$ is a derivation
of $B$ preserving $\fe$, and
\[
\ad_b\,\varepsilon_j(u_1\cdots u_j)=
\varepsilon_j\Bigl(\sum_iu_1\cdots[b,u_i]\cdots u_j\Bigr),
\]
so each $\varepsilon_j$ is $\ad(\fe)$-equivariant for the derivation
action on $S^j\fe$.

\begin{lemma}[symmetrization decomposition]\label{lem:symm}
\begin{enumerate}[label=\textup{(\arabic*)}]
\item $B=K\oplus\fe\oplus\bigoplus_{j=2}^{N-1}\varepsilon_j(S^j\fe)$;
every summand is $\ad(\fe)$-stable, and $\I(B)^m$ is the direct sum of
its intersections with the summands, for every $m$.
\item The image $J_B$ of $\I(KV^\circ)\I(KE^\circ)$ in $B$ equals
$\fv\cdot\I(B)$ and is $\bigoplus_{j=2}^{N-1}\varepsilon_j(\fv\cdot
S^{j-1}\fe)$, a direct sum of $\ad(\fe)$-stable subspaces of the summands
with $j\ge2$.
\item Consequently $s':=B/J_B=s(E^\circ,V^\circ)/\I^N$ decomposes as
$K\oplus\fe\oplus C$ with $C$ $\ad(\fe)$-stable and $\fe$ intact; the
projection $p_1:s'\to\fe$ along $K\oplus C$ satisfies
$p_1([u,b])=[p_1(u),b]$ for $u\in s'$, $b\in\fe$, and
$p_1(\I(s')^m)\subseteq\gamma_m(\fe)$; and the image of $\fv$ in $s'$ is
$\{V^\circ-1\}$, with $\log v=v-1$ for $v\in V^\circ$.
\end{enumerate}
\end{lemma}

\begin{proof}
(1) Filter $S^j\fe$ by the total degree induced from the lower central
filtration of $\fe$, and let $\varepsilon=\sum_j\varepsilon_j$ on
$K\oplus\fe\oplus\bigoplus_{j\ge2}S^j\fe$; it is filtered, and by
Lemma~\ref{lem:lazard} the associated graded of $\fe$ is $L:=L(E^\circ)$
inside $\gr B=u(L)$, which lives in degrees $\le N-1$, with trivial
$p$-map since $E^\circ$ has exponent $p$. So $\gr\varepsilon$ is the
symmetrization map $K\oplus\bigoplus_{j\le N-1}S^jL\to u(L)$. By the
Poincar\'e--Birkhoff--Witt theorem for restricted enveloping algebras
\cite[Chapter~V, \S7]{Jac62}, the associated graded of $u(L)$ for the length
filtration is $S(L)/(x^p:x\in L)$, and symmetrization $S^jL\to
U^{(j)}/U^{(j-1)}$ is the natural map onto its degree-$j$ part, which is
injective for $j\le p-1$; a triangular argument shows that
$\bigoplus_{j\le N-1}S^jL\to u(L)$ is injective. It is surjective in
degrees $\le N-1$, since a monomial of length $\ell$ in a homogeneous
basis of $L$ has degree $\ge\ell$ and is, modulo shorter monomials, the
symmetrization of the corresponding commutative monomial. So
$\gr\varepsilon$ is bijective, $\varepsilon$ is a strict isomorphism of
filtered spaces, and (1) follows.

(2) For $v\in V^\circ$, $v-1=(\log v)w_v$ with $w_v$ a unit, so
$J_B=\fv\cdot\I(B)=\sum_j\fv\cdot\varepsilon_j(S^j\fe)$. For $z\in\fv$
central and $j+1\le N-1$, inserting $z$ in each of the $j+1$ positions
gives $z\,\varepsilon_j(u_1\cdots u_j)=\varepsilon_{j+1}(zu_1\cdots u_j)$;
and $z\,\varepsilon_{N-1}(S^{N-1}\fe)\subseteq\I(B)^{N}=0$ since
$z\in\I(B)^2$. Stability holds because $\fv$ is an ideal of $\fe$ and
$\varepsilon_j$ is equivariant.

(3) The decomposition descends because $J_B$ is a direct sum of
subspaces of the summands with $j\ge2$; equivariance of $p_1$ is the
stability of the summands, and $p_1(\I(s')^m)\subseteq\gamma_m(\fe)$
follows from the strictness in (1) and $\fe\cap\I(B)^m=\gamma_m(\fe)$.
In $s'$ one has $(v-1)\I(s')=0$ by \eqref{eq:trunc0}, so $\log v=
(v-1)u_v=v-1$.
\end{proof}

\subsection{Rigidity below the ceiling}

\begin{theorem}[truncated relative rigidity]\label{thm:trunc}
Let $p$ be odd and $H$ of exponent $p$ and class $c\le p-2$. Then $H$ has
$(\RR_n)$ for every $n$ with $c+1\le n\le p-1$: for every
$\Gamma\in\UAut_n(H)$ there is $\beta\in\Aut(H)$, trivial on $H/\Phi(H)$,
with $\rho_n(\Gamma)=\rho(\beta)|_{\Vt}$.
\end{theorem}

\begin{proof}
Put $N=n+1\le p$. By Proposition~\ref{prop:E}(3), $W_{n+1}=\mho_1(E)$, so
$\Vt=V/\mho_1(E)=V^\circ\cong C(V)$. The kernel of $s(E,V)\to
s(E^\circ,V^\circ)$ is $\{\mho_1(E)-1\}$, which lies in $\I(s)^p\subseteq
\I(s)^N$ because $x^p-1=(x-1)^p$; hence $\st=s(E,V)/\I^N=
s(E^\circ,V^\circ)/\I^N=s'$ in the notation of Lemma~\ref{lem:symm}, and
$\widetilde{\mathcal V}=\{V^\circ-1\}$ is the image of $\fv$.

Let $\Gamma\in\UAut_n(H)$ and $\alpha=\tilde\Gamma\in\Aut(s')$, so that
$(\alpha-1)\I\subseteq\I^2$ by Lemma~\ref{lem:rho}(1). By
Lemma~\ref{lem:log}, applicable since $\I(s')^N=0$ and $N\le p$,
$D:=\log\alpha$ is a derivation of $s'$ with $D(\I^m)\subseteq\I^{m+1}$
and $\alpha=\exp D$; being a polynomial in $\alpha-1$, it preserves
$\widetilde{\mathcal V}=\fv$. Put $\tilde\mu:=p_1\circ D|_{\fe}:\fe\to\fe$.
By the equivariance of $p_1$,
\[
\tilde\mu[a,b]=p_1\bigl([Da,b]+[a,Db]\bigr)=[\tilde\mu a,b]+[a,\tilde\mu b],
\]
so $\tilde\mu$ is a derivation of $\fe$; it satisfies
$\tilde\mu(\gamma_m\fe)\subseteq p_1(\I^{m+1})\subseteq\gamma_{m+1}\fe$,
and $\tilde\mu|_{\fv}=D|_{\fv}$ because $D(\fv)\subseteq\fv\subseteq\fe$
and $p_1$ is the identity on $\fe$. Let $\mu$ be the derivation of
$\fe/\fv=\fl(H)$ induced by $\tilde\mu$; then $\mu(\fl(H))\subseteq
\gamma_2(\fl(H))$. Since $\fe$ has class at most $c+1\le n\le p-1$, we
have $\tilde\mu^{\,n}=0$, and $\Lambda:=\exp(\tilde\mu)=\sum_{k\le n-1}
\tilde\mu^k/k!$ is defined, with $(n-1)!$ invertible. It is a Lie
automorphism of $\fe$: expanding $\tilde\mu^k[a,b]=\sum_{i+j=k}\binom ki
[\tilde\mu^ia,\tilde\mu^jb]$ and using $[\gamma_{i+1}\fe,\gamma_{j+1}\fe]
=0$ for $i+j\ge n-1$ gives $\Lambda[a,b]=[\Lambda a,\Lambda b]$, and
$\Lambda$ is unipotent. It preserves $\fv$ and induces $\exp(\mu)$ on
$\fl(H)$, which is the identity modulo $\gamma_2$. By
Lemma~\ref{lem:lazard}, $\psi:=\exp\circ\Lambda\circ\log$ is an
automorphism of $E^\circ$ preserving $V^\circ$ and inducing on $H$ the
automorphism $\beta$ corresponding to $\exp(\mu)$, which is trivial on
$H/\gamma_2(H)=H/\Phi(H)$.

We compare the two actions on $V^\circ$. For $v\in V^\circ$, in $s'$ we
have $\log v=v-1$, and $\Lambda(v-1)=\exp(\tilde\mu)(v-1)=
\exp(D)(v-1)=\alpha(v-1)$, since $\tilde\mu|_{\fv}=D|_{\fv}$ and both
exponentials are the same truncated sum on $\fv$. Hence
$\log\psi(v)=\Lambda(\log v)$ has image $\alpha(v-1)=\alpha(v)-1$ in $s'$,
that is, $\psi(v)-1=\alpha(v-1)$ in $\widetilde{\mathcal V}$: under
$\Vt=V^\circ\cong\widetilde{\mathcal V}$, $\psi|_{V^\circ}=\rho_n(\Gamma)$.
On the other hand $\psi$ and $\beta^\circ$ of Lemma~\ref{lem:liftsE} are
automorphisms of $E^\circ$ preserving $V^\circ$ and inducing $\beta$, so
they agree on $V^\circ$ by Lemma~\ref{lem:welldef}, and $\beta^\circ$ acts
on $V^\circ=\Vt$ as $\rho(\beta)$ does. Thus $\rho_n(\Gamma)=
\rho(\beta)|_{\Vt}$.
\end{proof}

\begin{remark}\label{rem:calib}
The truncation degree $n+1$ is the smallest at which one can cut without
losing any of $\fe$, whose class is at most $c+1\le n$; what the
truncation destroys is exactly the power sector $\mho_1(E)=W_{n+1}$ and
the higher symmetric summands. The two inequalities used, $n+1\le p$ for
the logarithm and $\cl(\fe)\le p-1$ for the Lazard correspondence, are
both the hypothesis $c\le p-2$.
\end{remark}

\subsection{Rigidity at the ceiling}
At level $n=p$ the truncation no longer kills the power sector:
$W_{p+1}=1$ and $\Vt=V=C(V)\oplus\mho_1(E)$. The following lemma, whose argument is that of \cite[\S3]{UA}, shows that the
power sector is inert.

\begin{lemma}[shadow]\label{lem:shadow}
Let $\Gamma\in\UAut_p(H)$. Then $\rho_p(\Gamma)$ fixes $\mho_1(E)$
pointwise, and $\rho_p(\Gamma)(v)v^{-1}\in C(V)$ for every $v\in V$.
\end{lemma}

\begin{proof}
For $u\in\mho_1(E)\subseteq W_p$, Lemma~\ref{lem:rho}(1) gives
$\tilde\Gamma(u-1)-(u-1)\in\I(\st)^{p+1}=0$.

Put $\bar E=E/[E,E]\cong(C_{p^2})^d$, $\bar V=\varphi(V)$, and
$\bar s=s(\bar E,\bar V)/\I^{p+1}$. The map $KE\to K\bar E$ induces
$\psi:\st\to\bar s$, and $\Gamma$, being an automorphism of $A_{p+1}$,
preserves the ideal generated by the $ab-ba$, whose quotient is the
truncated group algebra of $F/[F,F]$; so $\Gamma$ induces an automorphism
$\bar\Gamma$ of $\bar s$ with $\psi\tilde\Gamma=\bar\Gamma\psi$. Write
$\Gamma(X_i)=X_i+\rho_i$ with $\rho_i\in\I^2$, so that
$\tilde\Gamma(\hat x_i)=\hat x_i(1+\hat x_i^{-1}\hat\rho_i)$ for the
images $\hat x_i$ of the $x_i$, and let $\sigma_i=\psi(\hat x_i^{-1}
\hat\rho_i)\in\I(\bar s)^2$. Let $v\in V$ be the class of a relator
$r\in R'$, and let $v'=\rho_p(\Gamma)(v)$, so that $\tilde\Gamma(\hat v)=
\hat v'$ by definition of $\rho_p$ (using $W_{p+1}=1$). Evaluating the
word $r$ in the commutative algebra $\bar s$,
\[
\psi(\hat v')=\bar\Gamma\psi(\hat v)=r\bigl(\psi(\hat x_i)(1+\sigma_i)\bigr)
=\psi(\hat v)\prod_i(1+\sigma_i)^{e_i(r)},
\]
where $e_i(r)$ is the exponent sum of $r$ in $x_i$. As $r\in R'\subseteq
F^{\,p}[F,F]$, $e_i(r)=pf_i$ with $f_i\in\Z$, and $(1+\sigma)^p=1+\sigma^p$
in $\bar s$, with $\sigma_i^{\,p}\in\I(\bar s)^{2p}=0$. Hence
$\psi(\hat v')=\psi(\hat v)$, that is $\varphi(v'v^{-1})-1=0$ in $\bar s$.
By Lemma~\ref{lem:PS} for the pair $(\bar E,\bar V)$, this means
$\varphi(v'v^{-1})\in\bar V\cap D_{p+1}(\bar E)=\mho_2(\bar E)=1$, so
$v'v^{-1}\in\ker(\varphi|_V)=C(V)$.
\end{proof}

\begin{theorem}[rigidity at the ceiling]\label{thm:ceiling}
Let $p$ be odd and $H$ of exponent $p$ and class $c\le p-2$. Then $H$ has
$(\RR_p)$.
\end{theorem}

\begin{proof}
Let $\Gamma\in\UAut_p(H)$ and let $\bar\Gamma\in\UAut_{p-1}(H)$ be its
image in $\Aut(A_p)$. The quotient map $\st_p\to\st_{p-1}=\st_p/\I^p$
carries $\{V-1\}$ onto $\{V-1\}/\{W_p-1\}$ and intertwines the two
actions, so $\rho_{p-1}(\bar\Gamma)$ is the map induced by $T:=
\rho_p(\Gamma)$ on $V/W_p$. Theorem~\ref{thm:trunc} at level $p-1$
(legitimate since $c+1\le p-1$) gives $\beta\in\Aut(H)$, trivial on
$H/\Phi(H)$, with $S:=\rho(\beta)$ inducing $\rho_{p-1}(\bar\Gamma)$ on
$V/W_p$. Now $V=C(V)\oplus\mho_1(E)$ with $W_p=\mho_1(E)$
(Proposition~\ref{prop:E}(3)). Both $T$ and $S$ preserve $C(V)$ ---
$T$ by Lemma~\ref{lem:shadow}, $S$ by Lemma~\ref{lem:liftsE} --- and both
fix $\mho_1(E)$ pointwise, by Lemma~\ref{lem:shadow} and
Lemma~\ref{lem:liftsE}. For $x\in C(V)$, $T(x)-S(x)$ lies in $C(V)$ and in
$W_p=\mho_1(E)$, hence is $0$. So $T=S$ on $V$.
\end{proof}

\section{The main theorems}\label{sec:main}

\begin{theorem}[the tower]\label{thm:tower}
Let $p$ be any prime and $G$ a finite $p$-group with $\nu=\nu(G)$. For
$2\le k\le\nu$ put $H_k=G/D_k(G)$. If $H_k$ has $(\RR_k)$ whenever
$D_k(G)\ne D_{k+1}(G)$, then $G$ is $(\nu+1)$-determined and has
$(\Lev_\nu)$; in particular $G$ is determined by $KG$ and
$\Aut(KG)=\Aut(G)\cdot\UAut(KG)$.
\end{theorem}

\begin{proof}
Fix a minimal presentation $F\to G$ and put $G_k=G/D_{k+1}(G)=F/R_k$ for
$k\ge1$, so that $R_k\le\Phi(F)$, $M_{k+1}(F)\subseteq R_k$, $d(G_k)=d$,
and $G_\nu=G$. We show by induction on $k$ that $G_k$ is
$(k+1)$-determined and has $(\Lev_k)$. For $k=1$, $G_1=G/\Phi(G)$ is
elementary abelian of rank $d$: $KX/\I^2$ determines $d(X)$, and the image
of $\Aut(G_1)$ in $GL_d(K)$ is everything. Let $k\ge2$. If
$D_k(G)=D_{k+1}(G)$ then $G_k=G_{k-1}$ and the claims follow from
\S\ref{subsec:det} and Lemma~\ref{lem:monotone}(1). Otherwise $G_k$ has
top Jennings layer $D_k(G_k)=D_k(G)/D_{k+1}(G)\ne1$ and base
$G_k/D_k(G_k)=G/D_k(G)=G_{k-1}=H_k$, which is $k$-determined, hence
$(k+1)$-determined, and has $(\Lev_{k-1})$ by induction, and $(\RR_k)$
by hypothesis. Lemma~\ref{lem:split}, Theorem~\ref{thm:descent} and
Theorem~\ref{thm:levi} conclude. The last sentence is \S\ref{subsec:det}
and Lemma~\ref{lem:monotone}(2).
\end{proof}

\begin{proof}[Proof of Corollary~\ref{cor:belowceiling}]
Put $\overline G=G/D_p(G)$. By Jennings' formula~\eqref{eq:jennings},
$D_n(\overline G)=D_n(G)/D_p(G)$ for $n\le p$, so $D_p(\overline G)=1$ and
$\nu(\overline G)\le p-1$. For $2\le k\le\nu(\overline G)$ the base
$\overline G/D_k(\overline G)=G/D_k(G)$ has exponent $p$, since
$\mho_1(G)\le D_p(G)\le D_k(G)$, and class at most $k-1\le p-2$, since
$\gamma_k(G)\le D_k(G)$. Hence it has $(\RR_k)$ by Theorem~\ref{thm:trunc},
and Theorem~\ref{thm:tower} shows that $\overline G$ is
$(\nu(\overline G)+1)$-determined, hence $p$-determined
(\S\ref{subsec:det}). Since $a-1\in\I^p$ for $a\in D_p(G)$, one has
$\Fp G/\I^p=\Fp\overline G/\I^p$, and likewise for $\Gt$. Thus
$\Fp\Gt/\I^p\cong\Fp\overline G/\I^p$, and $p$-determinacy of $\overline G$
gives $\Gt/D_p(\Gt)\cong\overline G$. The last sentence follows because an
isomorphism $\Fp G\cong\Fp\Gt$ maps $\I^p$ onto $\I^p$.
\end{proof}

\begin{proof}[Proof of Theorem~\ref{thmA}]
Let $p$ be odd, $D_{p+1}(G)=1$, so $\nu\le p$, and let $2\le k\le\nu$.
The base $H_k=G/D_k(G)$ has class at most $k-1$, since
$\gamma_k(G)\le D_k(G)$, and exponent $p$, since $\mho_1(G)\le D_p(G)\le
D_k(G)$ for $k\le p$. For $k\le p-1$ its class is at most $p-2$ and
$\cl(H_k)+1\le k\le p-1$, so $H_k$ has $(\RR_k)$ by
Theorem~\ref{thm:trunc}. For $k=p$ the hypothesis is that
$H_p=G/D_p(G)$ has class at most $p-2$, and $H_p$ has $(\RR_p)$ by
Theorem~\ref{thm:ceiling}. Theorem~\ref{thm:tower} applies.
\end{proof}

\begin{proof}[Proof of Corollaries~\ref{cor:expp} and~\ref{cor:expp2}]
By \eqref{eq:jennings},
\[
D_p(G)=\gamma_p(G)\,\mho_1(G),\qquad
D_{p+1}(G)=\gamma_{p+1}(G)\,\mho_1(\gamma_2(G))\,\mho_2(G).
\]
If $G$ has exponent $p$ and
class $<p$ then $D_p(G)=1$ and Theorem~\ref{thmA} applies with $\nu<p$.
If $\cl(G)\le p-2$, $\exp(G)\mid p^2$ and $\gamma_2(G)$ has exponent $p$,
then $D_{p+1}(G)=1$ and $D_p(G)=\mho_1(G)$, so $G/D_p(G)$ has class at
most $p-2$.
\end{proof}

\begin{proof}[Proof of Theorem~\ref{thmC}]
If $D_{p+1}(G)\ne1$, then by the formula above one of (I), (II), (III)
holds, and $\gamma_{p+1}\le D_{p+1}$, $\mho_1(\gamma_2)\le D_{2p}$,
$\mho_2(G)\le D_{p^2}$ give the bounds on $\nu(G)$. If $D_{p+1}(G)=1$,
then $G$ is not covered by Theorem~\ref{thmA} only if $D_p(G)\ne1$ and
$G/D_p(G)$, which has class at most $p-1$, has class exactly $p-1$; then
$\nu(G)=p$ and $\cl(G)\ge p-1$, and if $G$ has exponent $p$ then
$D_p(G)=\gamma_p(G)\ne1$, so $\cl(G)=p$. The isomorphism $G/D_p(G)\cong\Gt/D_p(\Gt)$ is
Corollary~\ref{cor:belowceiling}.
\end{proof}

The descent theorem also says where any failure must sit, in every
characteristic. Call $G$ \emph{Jennings determined} if it is
$(\nu(G)+1)$-determined.

\begin{proposition}[minimal failures]\label{prop:min}
Let $p$ be any prime, and suppose that some finite $p$-group is not
Jennings determined. Let $G$ be one of minimal order, $n=\nu(G)$,
$H=G/D_n(G)$. Then $n\ge2$, $H$ is $(n+1)$-determined, and hypothesis
\textup{(b)} of Theorem~\ref{thm:descent} fails for $H$ at level $n$:
either $H$ does not have $(\Lev_{n-1})$, or some unipotent displacement
$\rho_n(\Gamma)$, $\Gamma\in\UAut_n(H)$, is not realized by an automorphism
of $H$. In the second case the relation module $V/W_{n+1}$ of $H$
occupies at least two layers of the depth flag, and $H$ is not
relatively free.
\end{proposition}

\begin{proof}
Elementary abelian groups are $2$-determined, so $n\ge2$. By minimality
$H$ is Jennings determined; as $\nu(H)\le n-1$, it is $(n+1)$-determined,
which is hypothesis (a). The presentation $G=F/R$ satisfies
\eqref{eq:standing}, so Theorem~\ref{thm:descent} applies, and its
conclusion fails; hence (b) fails, and Lemma~\ref{lem:split} splits (b)
as stated. If the relation module lay in a single layer, the argument of
Corollary~\ref{cor:rohl} would give $(\RR_n)$ with $\beta=1$; and for $H$
relatively free both $(\Lev_{n-1})$ and the single-layer property hold.
\end{proof}

\section{Examples and scope}\label{sec:examples}

\subsection{The first groups outside previous results}\label{subsec:new}
Let $p\ge5$, let $F$ be free on $x,y$, and write $\ell$ for the free Lie
algebra on $x,y$ over $K$, $w=[x,y]$, $h_1=[w,x]$, $h_2=[w,y]$,
$z=[h_1,x]$; note $[h_2,x]=[h_1,y]$ by the Jacobi identity. The two
quotients of $\ell/\ell_{\ge5}$ by the ideals generated by $h_2$ and by
$h_2+z$ are the two filiform Lie algebras of dimension $5$, the first
graded and the second not: the centralizer of the derived algebra has
dimension $4$ in the first and $3$ in the second. Under the Lazard
correspondence they give two non-isomorphic groups $G_0,G_1$ of order
$p^5$, exponent $p$ and class $4$, with $L(G_0)\cong L(G_1)$ as graded
Lie algebras and $\gr(KG_0)\cong\gr(KG_1)$: the graded data do not
separate them, but Corollary~\ref{cor:expp} does, in agreement with
\cite{SS96}. Their common base $G_i/\gamma_4(G_i)$ is the group of
maximal class and order $p^4$, whose relation module occupies two layers
of the depth flag; so this pair is the smallest instance in which the
relation ideal does not lie in a single Jennings layer, and it is outside
Corollary~\ref{cor:rohl}. For $p=5$ the two groups are
$\mathrm{SmallGroup}(3125,30)$ and $(3125,33)$ and the base is
$\mathrm{SmallGroup}(625,7)$.

For $p=7$, an enumeration by the $p$-quotient algorithm \cite{ANUPQ} in
GAP \cite{GAP} finds exactly $32$ two-generator groups of exponent $7$,
class $5$ and order at most $7^7$: five of order $7^6$, all with
$G/\gamma_5(G)\in\{G_0,G_1\}$, and $27$ of order $7^7$. The relatively
free group of class $4$ and rank $2$ has order $7^8$, so none of them
has its relation subgroup in a single Jennings layer, and
Corollary~\ref{cor:rohl} does not apply; they have $D_4(G)\ne1$, order at
least $p^6$, class $5$ and $[G:Z(G)]\ge7^4$, and none is metacyclic.
Corollary~\ref{cor:expp} applies to all of them. To the author's knowledge
no previous result does.

\subsection{Consistency at $p=2$}\label{subsec:glmdr}
Theorem~\ref{thm:descent} holds for $p=2$, so it must be consistent with
the pair $G=\mathrm{SmallGroup}(512,456)$,
$\Gt=\mathrm{SmallGroup}(512,453)$ of \cite{GLMdR22}, and the
consistency is informative. One computes $\nu(G)=\nu(\Gt)=8$, $|D_8(G)|=|D_8(\Gt)|=2$, and
$G/D_8(G)\cong\Gt/D_8(\Gt)\cong\mathrm{SmallGroup}(256,3)=:H$, a group
with $\nu(H)=4$ and Jennings layer dimensions $(2,3,0,3)$. Since the two
bases are isomorphic, Step~1 of the proof of Theorem~\ref{thm:descent}
is not needed, and Steps~2--4, applied to an isomorphism
$KG\cong K\Gt$, show that if hypothesis (b) held for $H$ at level
$n=8$, then $G\cong\Gt$. Hence (b) fails for $H$ at $n=8$; in the notation of
Proposition~\ref{prop:count}, $[\mathcal G:\mathcal A]\ge2$.
The relation module of $H$ has dimension $5$ with depth flag of
dimensions $(5,5,5,3,3,3,3,3,0)$: two layers, in depths $3$ and $8$, and
the algebra fuses two $\Aut(H)$-orbits of hyperplanes, of sizes $2$ and
$1$, into one. This is the mechanism of \cite[Example~6.3]{HS07} in its
smallest form. In the terms of Theorem~\ref{thmC}, $G$ has exponent
$16\ge2^3$ and $\nu(G)=8\ge p^2$: the known counterexample lies in the
far case (III), and not one step beyond the range of this paper.

\subsection{Scope}\label{subsec:scope}
Theorem~\ref{thmA} does not touch class $\ge p+1$, derived subgroups of
exponent $p^2$, or elements of order $p^3$; in these directions
R\"ohl's theorem, which assumes nothing about the exponent, remains
strictly more general, and Theorem~\ref{thmC} measures how far each
lies beyond the ceiling. For $p=3$ the theorem covers $\nu(G)\le2$ and
the groups with $\nu(G)=3$ whose base $G/D_3(G)$ is abelian; Hertweck's
theorem covers all groups with $\nu(G)\le3$ and is the stronger one.
Everything here is over the prime field; the modular isomorphism
problem over all fields of characteristic $p$ is discussed in
\cite{MS25}, where even the case $D_3(G)=1$ is recorded as open.

The results of Section~\ref{sec:examples} were checked in GAP
\cite{GAP} with the packages ANUPQ~\cite{ANUPQ} and
AutPGrp~\cite{AutPGrp}, using the script provided as supplementary
material; none of the theorems depends on machine computation.

\section{Open problems}\label{sec:open}

\begin{problem}[the ceiling]
Extend Theorem~\ref{thmA} to bases $G/D_p(G)$ of class $p-1$, that is,
to case (IV) of Theorem~\ref{thmC}. The obstacle is precise: the
universal central Frattini extension $E$ of such a base has class $p$,
outside the Lazard correspondence, so Lemma~\ref{lem:lazard} is
unavailable for $E^\circ$. What the proof of Theorem~\ref{thm:trunc}
actually uses is $\gamma_p(E)=1$, which is implied by $\cl(H)\le p-2$
but is strictly weaker, and holds for instance when the relations of
$H$ in degree $p-1$ are dense enough for $[\ell_1,\bar N]=\ell_p$ in
the notation of \S\ref{subsec:new}.
\end{problem}

\begin{problem}[exponent $p^2$]
Extend the results to groups with $\mho_1(\gamma_2(G))\ne1$, where
$\nu(G)\ge2p$. Theorems~\ref{thm:descent} and~\ref{thm:levi} do not use
the exponent; what is missing is $(\RR_n)$ for bases of exponent $p^2$,
and by Remark~\ref{rem:calib} the logarithm cannot reach beyond
truncation degree $p$, so a different tool is needed.
\end{problem}

\begin{problem}[characteristic $2$]
Section~\ref{sec:descent} is valid for $p=2$. Determine, for bases $H$ of
order $2^8$, the index of $\rho(\Aut(H))|_{\Vt}$ in $\rho_n(\Aut_n(H))$;
by Proposition~\ref{prop:count} it bounds the number of groups sharing a
truncated modular group algebra with an extension of $H$, and its value
for $\mathrm{SmallGroup}(256,3)$ at $n=8$ is at least $2$.
\end{problem}

\section*{Acknowledgements}
The author completed his B.Sc.\ in Mathematics at the Instituto de
Matem\'atica e Estat\'istica, Universidade de S\~ao Paulo, in 2026, and
thanks the developers of GAP and of the packages ANUPQ and AutPGrp. No
funding was received for this work, and the author declares no
competing interests.

During the preparation of this work, the author consulted the large language
model Claude (Anthropic) for assistance in checking mathematical arguments,
exploring intermediate proof strategies, assisting in the development of the GAP
routines in Section~\ref{sec:examples}, and refining the exposition. All proofs,
computations, and arguments were independently analyzed, rigorously verified,
and written by the author, who assumes full responsibility for the mathematical
accuracy and content of the paper.

\section*{Data access statement}

The GAP code supporting the computations reported in
Section~\ref{sec:examples} is openly available in Zenodo at
\url{https://doi.org/10.5281/zenodo.22774195}.
No other data were generated.

\end{document}